\documentclass[12pt]{article}
\usepackage[utf8]{inputenc}
\usepackage{multirow}
\usepackage{theorem}%latexsym}
\usepackage{amsmath,amssymb}
\usepackage[T1]{fontenc}
\usepackage{graphicx}
\usepackage{color}
\usepackage[frenchb,english]{babel}
\usepackage{textcomp}
\usepackage{varioref}

\newcommand{\R}{\mathbb{R}}

\newcommand{\E}{\mathrm{E}}
\newcommand{\1}{1\hspace{-1.3mm}\mathrm{I}} 
\newcommand{\partiel}{\mathrm{d}}

 \theoremstyle{changebreak}
 \newtheorem{theorem}{Theorem}[section]

 \theoremstyle{changebreak}
 
 \theoremstyle{changebreak}
 
 \theoremstyle{changebreak}
 \newtheorem{lemma}[theorem]{Lemma}

 \theoremstyle{changebreak}
 \newtheorem{proposition}[theorem]{Proposition}

 \theoremstyle{changebreak}

 \theoremstyle{changebreak}
 \newtheorem{corollary}[theorem]{Corollary}

 \theoremstyle{changebreak}
 \newtheorem{remark}[theorem]{Remark}

 \author{Florent {\sc Nzissila}\footnote{nzissilaflorent@gmail.com }, 
	Octave {\sc Moutsinga}\footnote{\textbf{octave.moutsinga@univ-masuku.org}} and 
	Fulgence {\sc Eyi Obiang}\footnote{\textbf{feyiobiang@yahoo.fr}}
	}
 \date{ {\small  Universit\'e des Sciences et Techniques de Masuku\\
   Facult\'e des Sciences - Dpt Math\'ematiques et Informatique\\
  BP 943 Franceville, Gabon.\\
  URMI}}
	\title{Backwards semi-martingales into Burgers' turtulence}
\begin{document}
\bibliographystyle{plain}
%==================================================
\selectlanguage{english}
 \maketitle
\begin{abstract} 
In fluid dynamics governed by the one dimensional inviscid Burgers equation $\partial_t u+u\partial_x(u)=0$, the stirring is explained by the sticky particles model. A Markov process $([Z^1_t,Z^2_t],\,t\geq0)$ describes the motion of random turbulent intervals which evolve inside an other Markov process $([Z^3_t,Z^4_t],\,t\geq0)$, describing the motion of random clusters concerned with the turbulence. Then, the four velocity processes $(u(Z^i_t,t),\,t\geq0)$ are backward semi-martingales. 
\end{abstract}

\section{Introduction}
\indent
Burgers equation is a simplified version of Navier-Stocks equations in fluid mechanics. 
It is well known that the entropy solution of the one dimensional inviscid Burgers equation $\partial_t u+u\partial_x(u)=0$ can be interpreted as the velocity field of fluid particles which evolve following a sticky dynamics (\cite{Moutsinga-burgers-sticky1,Moutsinga-burgers-sticky2}). 
The aim of this paper is the study of fluid turbulence \textit{via} particular trajectories $s\mapsto y(s)$ of the fluid particles.

In the literature, the sticky particles dynamics was introduced, at a discrete level, by Zeldovich \cite{Zeldovich-Gravitational-instability} in order to explain the formation of large structures in the universe. That is a finite number of particles which move with constant velocities while they are not collided. All the shocks are inelastic following the conservation laws of mass and momentum. 

At a continuous level, the initial state of particles is given by the support of a non negative measure $\mu_0$. A particle starts from position $x$ with velocity $u_0(x)$ and mass $\mu_0(\{x\})$. The particles move with constant velocities and masses while not collided. All the shocks are inelastic, following the conservation laws of mass and momentum.

Now and in the rest of the paper, our purpose concerns only a motion on the real line. In their pioneering work, Sinai et al \cite{Generalized-variational-principles-E-Rykov-Sinai} made this construction when 
the particles are every where in $\R$, $u_0$ is continuous and the mass of any interval $[a,b]$ is computed with a positive density $f$, i.e. $\mu_0([a,b])=\int_a^bf(x)\partiel x$. At time $t$, a particle of position $x(t)$ has the mass $\mu_t(\{x(t)\})$ and the velocity $u_t(x(t))$, the momentum of any interval $[a(t),b(t)]$ is $\int_{a(t)}^{b(t)}u_t(x)\partiel\mu_t(x)$. The authors then solved the so called pressure-less gas system
$\partial_t \mu+\partial_x (u\mu)=0$, $\partial_t (u\mu)+\partial_x (u^2\mu)=0$. 

At the same time and independently, Brenier and Grenier \cite{sticky-part-Brenier-Grenier-2008} considered the case of particles confined in a in interval $[\alpha,\beta]$, i.e. $\mu_0([\alpha,\beta]^c)=0$. First, they remarked that at a discrete level of $n$ particles, the cumulative distribution function $M_n(a,t):=\mu_{n,t}(]-\infty,a])$ and the momentum $A_n(M_n(a,t)):=\int_{-\infty}^au_{n,t}(x)\partiel\mu_{n,t}(x)$ solve the so called scalar conservation law $\partial_t M_n+$ $\partial_x (A_n(M_{n}))=0$. Then, letting $n\to+\infty$, they got, at the continuous level, a limit $(M,A)$ solution of $\partial_t M+\partial_x (A(M))=0$. As a consequence, the Lebesgue-Stieltjes measure $\partial_x (A(M))$ is absolutely continuous w.r.t. $\partial_x M=:\mu_t$, of Randon-Nicodym derivative a function $u_t$, then $(\mu, u)$ solves again the above pressureless gas system.

In \cite{vpDermouneMoutsinga,  Convex-hullsMoutsinga}, Dermoune and Moutsinga constructed the sticky particles dynamics with initial mass distribution $\mu_0$, any probability measure, and initial velocity function $u_0$, any continuous and locally integrable function such that $u_0(x)=o(x)$ as $x\to\infty$. The authors united and generalized  previous works of \cite{Generalized-variational-principles-E-Rykov-Sinai,sticky-part-Brenier-Grenier-2008}. Moreover the particles paths define a Markov process $t\mapsto X_t$ solution of the ODE \begin{equation}\label{ODE}
\partiel X_t=u(X_t,t)\partiel t,\end{equation} and the velocity process $t\mapsto u(X_t,t)$ is a backward martingale.

In \cite{Convex-hullsMoutsinga, Moutsinga-burgers-sticky1}, Moutsinga extended the construction when $\mu_0$ is any non negative measure and $u_0$ has \textit{no positive jump}. He gave the description of different kinds of clusters $[\alpha(x,t),\beta(x,t)]$, \textit{i.e} the set of all the initial particles $y(0)$ which have the same position $y(t)=x$ at time $t$. The author showed in \cite{Moutsinga-burgers-sticky1} that if $\mu_0$ is the Lebesgue measure, then the velocity field $u$ is the entropy solution of the inviscid Burgers equation 
$\partial_t u+u\partial_x(u)=0$. In \cite{Moutsinga-burgers-sticky2}, Moutsinga showed the same connection  with Burgers equation when $u_0$ is non increasing and $\mu_0$ is the Stieltjes measure $-\partiel u_0$. In this case, the mass of any interval $[a(t),b(t)]$, at time $t$, is  
$u_t(a(t))-u_t(b(t))$ and its momentum is $(u_t^2(a(t))-u_t^2(b(t)))/2$.

The same year and independantly of Moutsinga, Eyink and Drivas \cite{spontaneous-Drivas-Eyink} considered $\mu_0=\lambda$, the Lebesgue measure (on a compact subset and renormalized to be a probability) and a random variable $\tau$ denoting the first shock time. They connected Burgers \textit{turbulence} to a Markov process  $s\mapsto Y_s$ solution of (\ref{ODE}). The authors also showed the anomalous dissipativity of $u$ along the turbulences' paths.

In fact, we show in section $\ref{section dpc}$ that this process is the path of some sticky particle.  Hence, the process of \cite{spontaneous-Drivas-Eyink} coincides with particular paths of the sticky particles process $X$ defined in \cite{Moutsinga-burgers-sticky1}.
A very interesting result of \cite{spontaneous-Drivas-Eyink} is that the process $s\mapsto u(Y(s),s)$ is a backward martingale, under the assumption of uniform distribution of $\tau$. Unfortunately, as it is shown in section $\ref{section turbulent traj}$, this assumption leads to the entire coincidence (undistinguishability) of both the processes $Y$ and $X$, so the martingale property of $s\mapsto u(Y(s),s)$ is obvious, since it was already stated in  \cite{Moutsinga-burgers-sticky1}. The construction of \cite{Moutsinga-burgers-sticky1} also allows us in subsection \ref{section dissipativity}, under more general assumptions than in \cite{spontaneous-Drivas-Eyink}, to show the anomalous dissipativity of the system governed by Burgers equation. 

Without the assumption of uniform distribution of $\tau$, the processes $Y$ and $X$ are in general distinguishable. We study this general case in section \ref{section turbulent traj} where we give the main results of this paper. The velocity function $u_0$ is not necessarily derivable nor even continuous as considered in \cite{spontaneous-Drivas-Eyink}, but it is allowed to have \textit{negative jumps}. We show that the process $t\mapsto u(Y_t,t)$ is no longer a backward martingale but a semi-martingale.  
Furthermore, we concentrate on the birth and evolution of turbulence. We define a \textit{turbulent interval} as a set $[\alpha,\beta]$ of initial positions of sticky particles from which rise a turbulence. The motion $s\mapsto [Z^1(s),Z^2(s)]$ of random turbulent interval is given by two backward Markov processes  $Z^1$ and $Z^2$ solutions of (\ref{ODE}). Moreover, the velocity processes 
$s\mapsto u(Z^1(s),s), u(Z^2(s),s)$ are semi-martingales.

First, we recall the definition and the main properties of the sticky particles model (\cite{Convex-hullsMoutsinga,Moutsinga-burgers-sticky1,Moutsinga-burgers-sticky2}).

\section{Flow and velocity field of sticky particles}\label{section dpc}
\subsection{The sticky particle dynamics}\label{the sticky dynamics}
The definition of one dimensional sticky particle dynamics requires a mass distribution $\mu$, any Radon measure (a measure finite on compact subsets) and a velocity function $u$, any real function such that the couple $(\mu,\,u)$ satisfies the \textit{Negative Jump Condition} (NJC) defined in $\cite{Convex-hullsMoutsinga}$. Precisely, consider the support  
$\mathcal{S}=\left\lbrace\,x\in\mathbb{R}\,:\, \mu(x-\varepsilon,\,x+\varepsilon)>0,\,\forall\varepsilon>0\right\rbrace$ of $\mu$ and the subsets  
$\mathcal{S}_{-}=\left\lbrace\,x\in\mathbb{R}\,:\, \mu(x-\varepsilon,\,x)>0\right\rbrace$, 
$\mathcal{S}_{+}=\left\lbrace\,x\in\mathbb{R}\,:\ \mu(x,\,x+\varepsilon)>0,\,\forall\varepsilon>0\right\rbrace$. 
Suppose that $u$ is $\mu$ locally integrable and consider the generalized limits $u^{-}$, $u^{+}$ :

\begin{align}\label{vitesse moins}
u^{-}(x)=\limsup\limits_{\varepsilon\rightarrow 0}\dfrac{\int_{[x-\varepsilon,\,x)}u(\eta)\mu(d\eta)}{\mu[x-\varepsilon,\,x)},\quad\forall\,x\in\mathcal{S}_{-},\\
\label{vitesse plus}
u^{+}(x)=\liminf\limits_{\varepsilon\rightarrow 0}\dfrac{\int_{(x,\,x+\varepsilon]}u(\eta)\mu(d\eta)}{\mu(x,\,x+\varepsilon]},\quad\forall\,x\in\mathcal{S}_{+}.
\end{align}
The \textit{Negative Jump Condition} requires that
\begin{equation}\label{njc}
u^{-}(x)\geq u(x)\,\,\forall\,x\in\mathcal{S}_{-},\quad u(x)\geq u^{+}(x)\,\,\forall\,x\in\mathcal{S}_{+}.
\end{equation}
In the whole paper, we mainly use $\mu_0=\lambda$, the Lebesgue measure. That's why we always suppose that the support $\mathcal{S}=\R$.
 
Considering particles of initial mass distribution $\mu_{0}$ and of initial velocity function $u_{0}$, their sticky dynamics is defined in $\cite{Moutsinga-burgers-sticky1}$, when the couple $(\mu_0,\,u_0)$ satisfies (\ref{njc}) and $x^{-1}u(x)\to0$ as $|x|\to+\infty$. The dynamics is characterized by a forward flow $(x,s,t)\mapsto\phi_{s,t}(x)$ defined on $\mathbb{R}\times\mathbb{R}_{+}\times\mathbb{R}_{+}$.
\begin{proposition}[Forward flow]\label{forward flow}
Suppose that $\mathcal{S}=\R$. For all $x,s,t$ : 
\begin{enumerate}
\item $\phi_{s,s}(x)=x$ and $\phi_{s,t}(\cdot)$ is non-decreasing and continuous.
\item The value $\phi_{s,t}(x)$ is the position after supplementary time $t$ of the particle which occupied the position $x$ at time $s$. More precisely :
\begin{equation}\label {flow-property}
\phi_{s,t}(\phi_{0,s}(y))=\phi_{0,s+t}(y)\;,\quad\forall\,y.
\end{equation}
\item If $\phi^{-1}_{0,t}(\left\lbrace x \right\rbrace )=:\left[ \alpha(x,0,t),\,\beta(x,0,t)\right]$ 
with $\alpha(x,0,t)<\beta(x,0,t)$, then
$$x=\frac{\int_{[ \alpha(x,0,t),\,\beta(x,0,t)]}(a+tu_0(a))\partiel\mu_0(a)}{\mu_0([ \alpha(x,0,t),\,\beta(x,0,t)])}\;.$$ 
In any case : 
$$x=\alpha(x,0,t)+tu_0(\alpha(x,0,t))=\beta(x,0,t)+tu_0(\beta(x,0,t))\;.$$
\item 
If $\mu_0([\alpha(x,0,t),\,y])>0$ and $\mu_0(]y,\,\beta(x,0,t)])>0$, then
\begin{align*}
\frac{\int_{]y, \beta(x,0,t)]}(a+tu_0(a))\partiel\mu_0(a)}{\mu_0(]y, \beta(x,0,t)])}
\leq x\leq
\frac{\int_{[ \alpha(x,0,t),y]}(a+tu_0(a))\partiel\mu_0(a)}{\mu_0([\alpha(x,0,t),y])}\;.
\end{align*}
\item If $s\leq t$, then 
\begin{align*}
\phi_{0,s}(\alpha(x,0,t))&=\alpha(x,0,t)+su_0(\alpha(x,0,t))\;,\\
\phi_{0,s}(\beta(x,0,t))&=\beta(x,0,t)+su_0(\beta(x,0,t))\;.
\end{align*}
\item For any compact subset $K=[a,b]\times[0,T]$, consider $A_T=\alpha(\phi_{s,T}(a),s,T)$, $B_T=\beta(\phi_{s,T}(b),s,T)$ and the probability $\mu_s^K=\frac{\1_{[A_T,B_T]}}{\mu_s([A_T,B_T])}\mu_s$. The sticky particle dynamics induced by $(\mu_s^K,u_s)$, during time interval $[0,T]$, is characterized by the restriction of the function $(y,t)\mapsto\phi_{s,t}(y)$ on $[A_T,B_T]\times[0,T]$. 
\end{enumerate} 
\end{proposition} 

The latter means that the restriction of flow on a compact subset of space-time does not depend of the whole matter, but only on the restriction of the matter (distribution) on a compact subset of space states.

\begin{figure}[h!]
\centering
\includegraphics[width=0.7\linewidth]{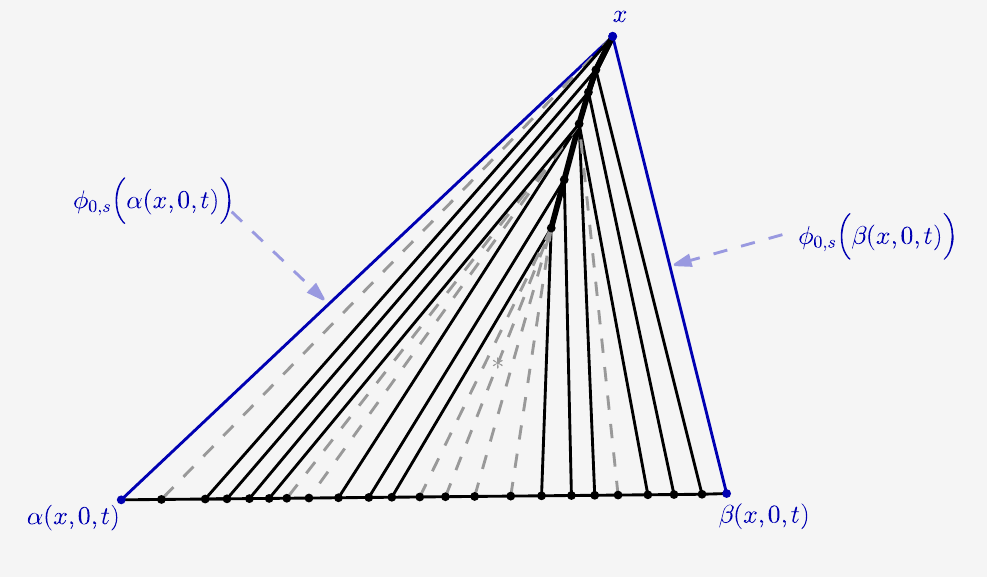}
\caption[delta-shock]{\scriptsize {The blue line on the left (resp right) of the middle shok wave represent the trajectory of the particle which start at the position $\alpha(x,0,t)$ (resp $\beta(x,0,t)$) which is the trajectory of $[0,t]\ni s\mapsto \phi_{0,s}(\alpha(x,0,t))$ (resp  $[0,t]\ni s\mapsto \phi_{0,s}(\beta(x,0,t))$)}}
\label{delta-shock}
\end{figure}

 Assertion 5 shows that the graphs $[0,t]\ni s\mapsto \phi_{0,s}(\alpha(x,0,t)),\,\phi_{0,s}(\beta(x,0,t))$ draw a delta-shock, well known in the literature (Figure \ref{delta-shock}).  \\
What about the velocity?

 \begin{proposition}[Flow derivative]\label{th derivative}
Suppose that $\mathcal{S}=\R$.
\begin{enumerate}
\item For all $y,s$, the function $t\mapsto\phi_{s,t}(y)$ has everywhere left hand derivatives and right hand derivatives.  Now and after, the notation
$\dfrac{\partial}{\partial t}\phi_{s,t}(y)$ stands for the right hand derivative.  

\item 
There exists a function
 $(x,t)\mapsto u_t(x)$ such that for all $(y,t)$: 
$$\dfrac{\partial}{\partial t}\phi_{0,t}(y)=u_t(\phi_{0,t}(y))\;.$$
\item For any compact subset $K=[a,b]\times[0,T]$, consider $A_T=\alpha(\phi_{s,T}(a),s,T)$, $B_T=\beta(\phi_{s,T}(a),s,T)$ and the probability $\mu_s^K=\frac{\1_{[A_T,B_T]}}{\mu_s([A_T,B_T])}\mu_s$. 

Using conditional expectation under $\mu_s^K$, we have
\begin{equation} \label{derivative of the flow}
\forall\,(y,t)\in K,\quad\dfrac{\partial}{\partial t}\phi_{s,t}(y)=\mathbb{E}_{\mu_s^K}\left[u_{s}|\phi_{s,t}(\cdot)= \phi_{s,t}(y)\right]. 
\end{equation}
\end{enumerate} 
\end{proposition} 

We call \textit{a cluster at time $t$} all interval
of the type $\left[ \alpha(x,0,t),\,\beta(x,0,t)\right]$. 
The last assertion of proposition \ref{forward flow} implies an important property on the velocity of a cluster.

\begin{corollary}\label{th velocity GVP}
Suppose that $\mathcal{S}=\R$. Let $(x,t)\in\R\times\R_+$.
\begin{enumerate}
\item If $\alpha(x,0,t)<\beta(x,0,t)$, then
$$u_t(x)=\frac{\int_{[ \alpha(x,0,t),\,\beta(x,0,t)]}u_0(a)\partiel\mu_0(a)}{\mu_0([ \alpha(x,0,t),\,\beta(x,0,t)])}\;.$$ 
If $\alpha(x,0,t)=\beta(x,0,t)$, then $u_t(x)=u_0(\alpha(x,0,t))$.
\item $u_0(\beta(x,0,t))\leq u_t(x)\leq u_0(\alpha(x,0,t))$.\\ 
If $\mu_0([\alpha(x,0,t),\,y])>0$ and $\mu_0(]y,\,\beta(x,0,t)])>0$, then
\begin{align*}
\frac{\int_{]y, \beta(x,0,t)]}u_0(a)\partiel\mu_0(a)}{\mu_0(]y, \beta(x,0,t)])}\leq u_t(x)\leq
\frac{\int_{[ \alpha(x,0,t),y]}u_0(a)\partiel\mu_0(a)}{\mu_0([\alpha(x,0,t),y])}\,.
\end{align*}
\item $u_t^-(x)=u_0(\alpha(x,0,t))$ and $u_t^+(x)=u_0(\beta(x,0,t))$).
\item If $u_0(\alpha(x,0,t))=u_t(x)$ or $u_0(\beta(x,0,t))=u_t(x)$, then  $\alpha(x,0,t)=\beta(x,0,t)$. 
\item For all $t\geq0$, we have $u_t(x)=o(x)$ as $|x|\to +\infty$. For all $t>0$ :
\begin{align*}
\underset{\underset{y<x}{y\to x}}{\lim}\,u_t(y)&=u_t^-(x)=u_0(\alpha(x,0,t))\;,\\
\underset{\underset{y>x}{y\to x}}{\lim}\,u_t(y)&=u_t^+(x)=u_0(\beta(x,0,t))\;.
\end{align*}

\end{enumerate} 
\end{corollary}

\subsection{Markov and martingale properties}\label{section Markov and martingale}
Let $(\mu_{0},\,u_{0})$ be as in theorem $\ref{forward flow}$. On abstract measure space $(\Omega,\,\mathcal{F},\,P)$ we define a measurable function $X_{0}\,:\;\Omega\longrightarrow\mathbb{R}$ with image-measure $P\circ X^{-1}_{0}=\mu_{0}$. In practice,  $(\Omega,\,\mathcal{F},\,P)=(\mathbb{R},\mathcal{B}(\mathbb{R}),\mu_{0})$ and $X_0$ is the identity function. For all  $t\geq0$, we set $X_{t}=\phi_{0,t}(X_{0})$. As a consequence of theorem $\ref{forward flow}$, we have the following :

\begin{proposition}[Markov and martingale property]\label{th martingale and markov prop}
\begin{enumerate}
\item $\forall s,t$, we have 
\begin{equation}\label{markov-property}X_{s+t}=\phi_{s,t}(X_{s})
\end{equation}
\item If $u_{0}$ is $\mu_0$ integrable, then under the measure $\mu_0$ (or $P$) :
\begin{equation}\label{velocity-process}
 \dfrac{\partiel }{\partiel t}X_{t}=\mathbb{E}[u_{0}(X_0)|X_{t}]=u_{t}(X_{t}).
\end{equation}
\item For any compact $K=[a,b]\times[0,t]$, consider $A_t=\alpha(\phi_{0,t}(a),0,t)$, $B_t=\beta(\phi_{0,t}(a),0,t)$ and the probability $\mu_0^K=\frac{\1_{[A_t,B_t]}}{\mu_0([A_t,B_t])}\mu_0$. \\
If $\phi_{0,t}(a)\leq X_{t}\leq\phi_{0,t}(b)$, then then under the conditional probability $\mu_0^K$ (or knowing $A_t\leq X_0\leq B_t$) :
\begin{equation}
 \dfrac{\partiel }{\partiel t}X_{t}=\mathbb{E}_{\mu_0^K}[u_{0}(X_0)|X_{t}]=u_{t}(X_{t}).
\end{equation}
\item If $u_{0}$ is $\mu_0$ integrable, then under the measure $\mu_0$ (or $P$) :
\begin{equation}\label{martingale-egality}
u_{t+s}(X_{t+s})=\mathbb{E}[u_{t}(X_{t})|\mathcal{F}_{t+s}]\;,\quad\mbox{with}\quad\mathcal{F}_{t}=\sigma(X_{u},u\geq t).
\end{equation}
\item For any compact $K=[a,b]\times[0,t+s]$, consider $A_{t+s}=\alpha(\phi_{0,t+s}(a),0,t+s)$, $B_{t+s}=\beta(\phi_{0,t+s}(a),0,t+s)$ and the probability $\mu_0^K=\frac{\1_{[A_{t+s},B_{t+s}]}}{\mu_0([A_{t+s},B_{t+s}])}\mu_0$. \\
If $\phi_{0,t+s}(a)\leq X_{t+s}\leq\phi_{0,t+s}(b)$, then then under the conditional probability $\mu_0^K$ (or knowing $A_{t+s}\leq X_0\leq B_{t+s}$) :
\begin{equation}
 \dfrac{\partiel }{\partiel t}X_{t+s}=\mathbb{E}_{\mu_0^K}[u_{s}(X_s)|\mathcal{F}_{t+s}]=u_{t+s}(X_{t+s}).
\end{equation}
\end{enumerate}
\end{proposition}

\begin{remark}
Contrary to the conjecture of \cite{spontaneous-Drivas-Eyink} (page $411$) the properties of proposition \ref{th martingale and markov prop} do not ensure $u$ to be the inviscid Burgers solution. 
\end{remark}

\begin{figure}[h!]
\centering
\includegraphics[width=0.9\linewidth]{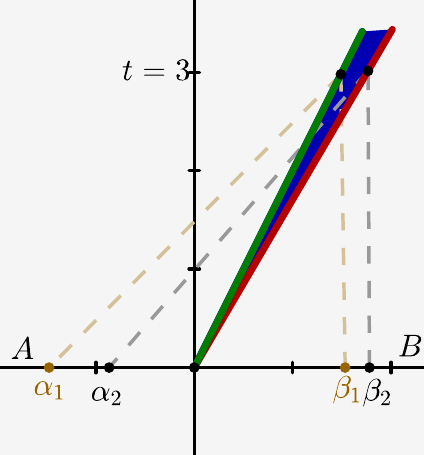}
\caption[cone-no-Burgers]{\scriptsize   {The green line is the discontinuous line in a lagrangian interval $[A,\,B]$ of the function $u$ described at example $1$. The red line is the discontinuous line in a Lagrangian interval $[A,\,B]$ of the function $u$ described at example $2$. $[\alpha_{1},\,\beta_{1}]$ is the cluster which contain the point zero for the discontinuous line in the first example and  $[\alpha_{2},\,\beta_{2}]$ is the cluster which contain the point zero for the discontinuous line in the second example. The two functions $u(x,t)$ described in the example $1$ and example $2$ don't coincide in the region of the plan which is red, blue and green. This region of the plan is not negligible for the Lebesgue measure. }}
\label{one process-different velocities}
\end{figure}

We can indeed give two examples (Figure \ref{one process-different velocities}) where the matter is initially confined in an interval $[-A,A]$ and 
$$u_0(x)=\left\{
\begin{array}{rll}
1 & \mbox{if} & -A\leq x\leq0\\
0 & \mbox{if} & 0<x\leq A.
\end{array}\right. 
$$ 
\begin{description}
\item{\textbf{Example 1 (dealing with Burgers equation) :}} $\mu_{0}=\lambda_{[-A,A]}$  (the Lebesgue measure on $[-A,A]$). We have a single discontinuity line (shock wave) $t\mapsto t/2$ starting at position 0, with velocity $1/2$. At time $t$, the position $x=t/2$ is the one of the cluster $[-t/2,t/2]$. Moreover,
$$
u(x,t)=\left\{
\begin{array}{rll}
1 & \mbox{if} & A+t\leq x<t/2\\
1/2 & \mbox{if} & x=t/2\\
0 & \mbox{if} & t/2<x\leq A.
\end{array}\right. 
$$

\item{\textbf{Example 2 (not dealing with Burgers equation) :}} $\mu_{0}=2\lambda_{[-A,0]}+\lambda_{[0,A]}$. We have a single discontinuity line (shock wave) $t\mapsto (2-\sqrt{2})t$ starting at position 0, with velocity $2-\sqrt{2}$. At time $t$, the position $x=(2-\sqrt{2})t$ is the one of the cluster $[(1-\sqrt{2})t,\,(2-\sqrt{2})t]$. Moreover,
$$
u(x,t)=\left\{
\begin{array}{rll}
1 & \mbox{if} & A+t\leq x<(2-\sqrt{2})t\\
(2-\sqrt{2}) & \mbox{if} & x=(2-\sqrt{2})t\\
0 & \mbox{if} & (2-\sqrt{2})t<x\leq A.
\end{array}\right. 
$$
\end{description} 

\subsection{Link with Burgers equation}\label{burger link}
The inviscid Burgers solution of initial data $u_{0}$ is connected to sticky particles in two well known cases ($\cite{Moutsinga-burgers-sticky1,Moutsinga-burgers-sticky2}$) as follows : 

\begin{proposition}\label{ lambda mu linkburgers}
\begin{enumerate}
\item When $\mu_{0}=\lambda$ (the Lebesgue measure) and $(\mu_{0},\,u_{0})$ satisfies the NJC, the function $u(\cdot,t)=u_{t}$ defined at the relation $(\ref{velocity-process})$ is the the entropy solution of the inviscid Burgers equation with initial data $u_{0}$. Furthermore the distribution of matter at time $t\geq 0$ is given by the relation 
\begin{equation}\label{t-distribution-lambda}
\mu_{t}=\lambda\circ X^{-1}_{t}=\lambda-t\partial_{x}u(\cdot,t)
\end{equation} where $\partial_{x}u(\cdot,t)$ is the Stieltjes measure over $u(\cdot,t)$.
\item When $\mu_{0}=-\partiel u_{0}$ (the Stieltjes measure over $u_{0}$) and $u_{0}$ is non-increasing,  the function $u(\cdot,t)=u_{t}$ defined at the relation $(\ref{velocity-process})$ is the the entropy solution of the inviscid Burgers equation with initial data $u_{0}$. Furthermore the distribution of matter at time $t\geq 0$ is given by the relation 
\begin{equation}\label{t-distribution-stieljes}
\mu_{t}=(-\partiel u_{0})\circ X^{-1}_{t}=-\partial_{x}u(\cdot,t)
\end{equation}
\end{enumerate}
\end{proposition}

\subsection{Dissipativity}\label{section dissipativity}
\begin{proposition}
Let $u$ be the velocity field of sticky particles and $\mu$ be the distribution field. 
Let $\psi$ be a convex function and $0\leq s<t$. 
\begin{enumerate}
\item If  
$\int|\psi(u(x,s))|\partiel\mu_s(x)+\int|\psi(u(x,t))|\partiel\mu_t(x)<+\infty$, then
\begin{align*}
\int\psi(u(x,t))\partiel\mu_t(x)\leq\int\psi(u(x,s))\partiel\mu_s(x)\;,
\end{align*}
which is equivalent to
\begin{align*}
\int\psi(u(\psi_{0,t}(a),t)\partiel \mu_0(a)\leq\int\psi(u(\psi_{0,s}(a),s)\partiel \mu_0(a)\;.
\end{align*}
\item Let $u$ be the entropy solution of the inviscid Burgers equation.\\
Consider $\mu_0=\lambda$, the Lebesgue measure and suppose that $\lambda$-essentially, $u_0(x)\to0$ as $x\to\pm\infty$.  
If
$\int[|\psi(u_0(x))|+|\psi(u(x,t))|]\partiel x<+\infty$, then
\begin{align*}
\int\psi(u(x,t))\partiel x=\int\psi(u(x,t))\partiel\mu_t(x)\leq\int\psi(u_0(x))\partiel x\;.
\end{align*}
\end{enumerate}
\end{proposition}

\textbf{Proof.} Let us use the probabilistic notations of subsection \ref{section Markov and martingale} : 
$(X_t,P)=(\psi_{0,t},\partiel\mu_0)$. 
For assertion 1), we first suppose that $u_0$ is $\mu_0$ integrable. We use Jensen inequality :
\begin{align*}
\psi(u(X_t,t))=
\psi(\E[u(X_s,s)|X_t])\leq \E[\psi(u(X_s,s))|X_t]\;.
\end{align*}
So 
\begin{align*}
\int\psi(u(x,t))\partiel\mu_t(x)=\E[\psi(u(X_t,t))]\leq\E[\psi(u(X_s,s))]=\int\psi(u(x,s))\partiel\mu_s(x)\;.
\end{align*}
 If $u_0$ is not $\mu_0$ integrable, we slightly modify the proof.
Consider $y<z$ and $a_s=\alpha(y,s,t-s)$,  $b_s=\beta(z,s,t-s)$ and $K=[a_0,b_0]\times[0,t]$. Using Jensen inequality, we have
\begin{align*}
\psi(u(X_t,t))\1_{[y,z]}(X_t)&=
\psi(\E_{\mu_0^K}[u(X_s,s)|X_t])\1_{[y,z]}(X_t)\\
&\leq \E_{\mu_0^K}[\psi(u(X_s,s))|X_t]\1_{[y,z]}(X_t)\;.
\end{align*}
So $\E_{\mu_0^K}[\psi(u(X_t,t))\1_{[y,z]}(X_t)]\leq\E_{\mu_0^K}[\psi(u(X_s,s))\1_{[y,z]}(X_t)]$.\\
Since $\1_{[y,z]}(X_t)=\1_{[a_s,b_s]}(X_s)$, this leads to 
\begin{align*}
\int_y^z\psi(u(x,t))\partiel\mu_t(x)\leq\int_{a_s}^{b_s}\psi(u(x,s))\partiel\mu_s(x)\;.
\end{align*}
When $y\to-\infty$ and $z\to+\infty$, we have $a_s\to-\infty$ and $b_s\to+\infty$. So we get
\begin{align*}
\int\psi(u(x,t))\partiel\mu_t(x)\leq\int\psi(u(x,s))\partiel\mu_s(x)\;.
\end{align*}

 Now, if $u$ is moreover the entropy solution of the inviscid Burgers equation and $\mu_0=\lambda$, 
then (\ref{t-distribution-lambda}) holds. Hence
\begin{align*}
\int\psi(u(x,t))\partiel x&=\int\psi(u(x,t))\partiel\mu_t(x)+t\int\psi(u(x,t))\partial_xu(x,t)\\
&=\int\psi(u(x,t))\partiel\mu_t(x)+t[\Psi(u(+\infty,t))-\Psi(u(-\infty,t))]\;,
\end{align*}
where $\Psi$ is any primitive of $\psi$. 
As $x\to\pm\infty$, the $\lambda$-essential limit of $u_0$ is 0. Thus $u(x,t)\to0$ as $x\to\pm\infty$.
The proof ends using assertion 1 with $s=0$ and $\mu_0=\lambda$.

\section{Dynamics of Burgers turbulence}\label{section turbulent traj}
In this section, we study the sticky particles dynamics from the point of view of turbulence. There emerge four more Markov processes solution of (\ref{ODE}). A remarkable property is that their velocities are backward semi-martingales.

We always suppose that the support $\mathcal{S}=\R$.

\subsection{Turbulence time and semi-martingales}\label{section Moutsinga-EyinkDrivas}

Following the preoccupation of \cite{spontaneous-Drivas-Eyink}, we consider the first shock time and related delta-shocks. 
Consider the left hand limit function $u^{-}(\cdot,t)$ and the right hand limit function $u^{+}(\cdot,t)$. Remark that if a particle of initial position $a$ is in a shock at time $t$ and position $x$, then $x=\phi_{0,t}(a)$ and $u^{-}(x,t)\neq u^{+}(x,t)$. That's why we define its first shock time as 
\begin{eqnarray}\label{tau def}
\tau(a)=\inf\left\lbrace t:\,u^{-}(\phi_{0,t}(a),t)\neq u^{+}(\phi_{0,t}(a),t)\right\rbrace .
\end{eqnarray}
In fact, as described in the sequel, this is more a turbulence time than a shock time. 
We define a turbulent interval containing $a$ as the greatest interval $[A,B]$ of initial positions of particles containing $a$ and which have first turbulence at same time and same position : $\tau(a')=\tau(a)=:t$ and $\phi_{0,t}(a')=\phi_{0,t}(a)$ $\forall\,a'\in[A,B]$. Because of the regularity of $u_0$ and $\phi$, a turbulent interval is effectively closed. 
 
\begin{figure}[h!]
\centering
\includegraphics[width=1.2\linewidth]{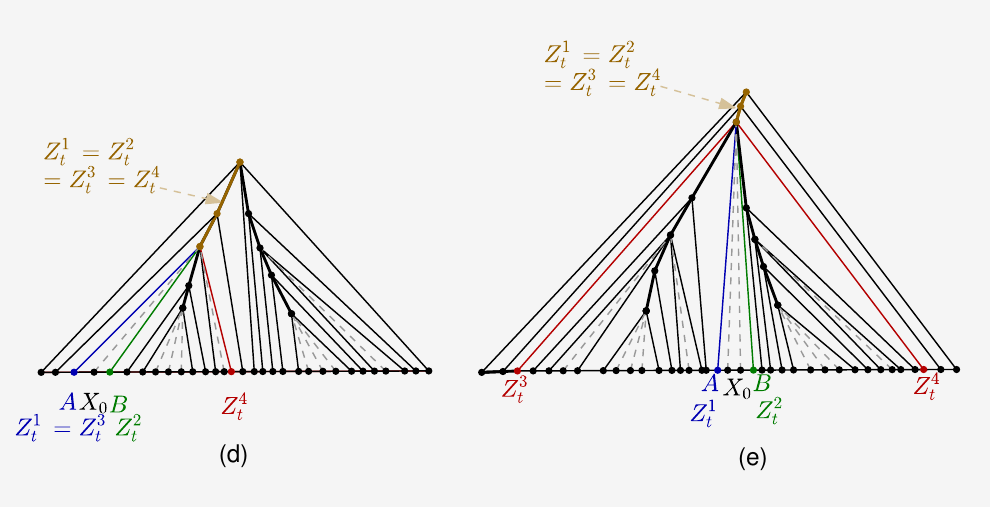}
\caption[two merger shocks]{ \scriptsize   {Merger of Two Shocks waves}}
\label{Shocks-waves}
\end{figure}

At time $\tau(a)$, there is a turbulence located at $x=\phi_{0,\tau(a)}(a)=a+\tau(a)u_0(a)$. The triangle, in the space-time representation, delimited by $A,B,x$ in known in the literature as a delta-shock (see Figure \ref{Z1,2t-trajectories}). As it is related to first shock, we call it a \textit{prime-delta-shock}. Thus, in the space-time representation, the turbulences are entirely conditioned by prime-delta-shocks. 
\begin{figure}[h!]
\centering
\includegraphics[width=1.\linewidth]{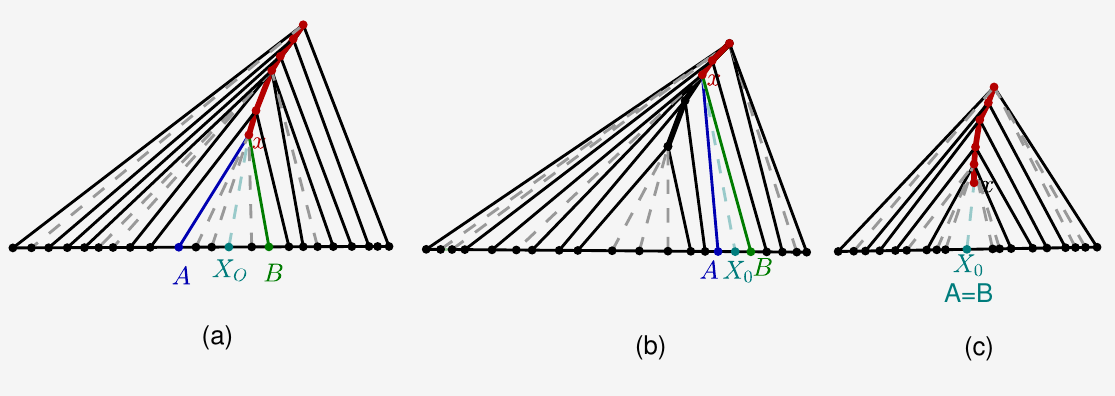}
\caption{{\scriptsize   {The blue curve and red curve in $(a)$ and $(b)$ (resp. the green curve and red curve in $(a)$ and $(b)$) represent the trajectory $Z^{1}_{t}$ (resp. $Z^{2}_{t}$) for $X_{0}\in [A,\,B]$. $x$ represent the position of the shock of particles from  $[A,\,B]$: the first shock position. In case (c), the red curve represents one trajectory of $Z^ {1}_{t}=Z^{2}_{t}$ in the discontinuous line. In this case, the point $x$ is a turbulent point. Immediately after the point $x$, there is shocks, but there is not a real cluster which get the position $x$ .}}} 
\label{Z1,2t-trajectories}
\end{figure}
Turbulent intervals bring suddenly positive masses to shocks. 
If the turbulent interval $[A,B]$ is a cluster at time $\tau(A)$, a turbulence of length $B-A$ rises from $[A,B]$; it is born at time $\tau(A)$ (see Figure \ref{Z1,2t-trajectories} a). If moreover $A=B$, then the turbulence is not detectable when it appears, an infinitesimal colliding (agglomeration) process starts at position $x=A+\tau(A)u_0(A)$ (see Figure \ref{Z1,2t-trajectories} c).
If the turbulent interval $[A,B]$ is not a cluster at time $\tau(A)$, it simply aggregates (in the same way as above) a turbulence which was born earlier from another turbulent interval (see Figure \ref{Z1,2t-trajectories} b) : 
$\exists a_0\not\in[A,B]$, $\exists t<\tau(A)$ such that 
$$\phi_{0,\tau(A)}(a_0)=\phi_{0,\tau(A)}(A)\; \mbox{ and }
\{a_0\}\subsetneq\{a': \phi_{0,t}(a')=\phi_{0,t}(a_0)\}.$$
Now, for any turbulent interval $[A,B]$, consider $x=A+\tau(A)u_0(A)$. We define four processes : 
if $X_0\in[A,B]$, then $\forall\,t\geq0$,
\begin{align}
Z^1_t=\phi_{0,t}(A)\;,\quad& Z^2_t=\phi_{0,t}(B)\;,\\
Z^3_t=\phi_{0,t}(\alpha(x,0,\tau(A))),\quad& Z^4_t=\phi_{0,t}(\beta(x,0,\tau(A)))\;.
\end{align}
See Figure \ref{Shocks-waves} for illustration.
It is clear that the set of all the turbulent intervals is a partition of the whole initial state of particles. Hence, the stochastic processes $Z^1,Z^2,Z^3$ and $Z^4$ are well defined. They all give new point of views of the sticky particle dynamics. The process $[Z^1,Z^2]$ describes the motion turbulent interval while $[Z^3,Z^4]$ describes the motion of all the particles which are concerned in the turbulence. 
Moreover, $[Z^1,Z^2]\subset[Z^3,Z^4]$. 

The measurability of $Z^1$ and $Z^2$ comes from the fact that $\sigma(Z^1_0,Z^2_0)\subset\sigma(X_0)$. Indeed,
for all $x\in\R$, the event $\{Z^1_0\leq x\}=\{X_0\leq B\}$ and $\{Z^2_0\leq x\}=\{X_0\leq A\}$, where $[A,B]$ is the turbulent interval such that $A\leq x\leq B$. For $Z^3$ and $Z^4$, we recall that 
for all $t$, the functions $\alpha(\cdot,0,t)$ and $\beta(\cdot,0,t)$ are non decreasing, so they are Borel functions; then 
$\alpha(X_t,0,t)$, $\beta(X_t,0,t)$ are $\sigma(X_t)$ measurable. Furthermore, all the paths of $t\mapsto\alpha(X_t,0,t)$, $\beta(X_t,0,t)$ are \textit{c\`adl\`ag}. The proof of the measurabilty of $Z^3_0$ and $Z^4_0$ is achieved by lemmas \ref{th measurability} and \ref{th stopping time}, since 
$$Z^3_0=\alpha(X_{\gamma},0,\gamma),\quad Z^4_0=\beta(X_{\gamma},0,\gamma),$$ 
where 
\begin{align}\label{gamma def}
\gamma=\tau(X_0)=\inf\left\lbrace t:\,u^{-}(X_t,t)\neq u^{+}(X_t,t)\right\rbrace 
\end{align}
is the first shock time of $X_0$. Remark that if $\Omega=\R$ and $X_0$ is the identity function, then $\tau=\gamma$. %
By construction, we have a.e :

\begin{proposition}[Random delta-shock] \label{th turbulent flows}
Suppose that the support $\mathcal{S}=\R$.
\begin{enumerate}
\item Let $Z$ stand independently for $Z^1,Z^2,Z^3$ or $Z^4$. 
\begin{align*}
\forall\,t,s\geq0,\quad Z_{s+t}=\phi_{s,t}(Z_{s})\;,\quad\frac{\partiel }{\partiel t}Z_{t}=u(Z_{t},t)\;.
\end{align*}
\item $\tau(Z_{0})=\tau(X_{0})=\gamma$ and $\forall\,t\geq\gamma$, $Z_{t}=X_t$;
\begin{align*}
\forall\,t\leq\gamma,\quad&Z_{t}=Z_{0}+t u_0(Z_{0})\;,\\
&Z^3_{t}\leq Z^1_{t}\leq X_t\leq Z^2_{t}\leq Z^4_{t}\;.
\end{align*}
\end{enumerate}
\end{proposition}
The graphs $[0,\gamma]\ni t\mapsto\,Z^1,Z^2,Z^3,Z^4$ then draw nested delta-shocks 
(see figure \ref{Shocks-waves}). 

Now we recall some properties well known in the the theory of stochastic processes (\cite{Ashkan-theory-of-processses,Xi-Geng-stochastic-calculs}) for forward martingales and non-decreasing filtrations. By inversion of time, the following holds.  
 
\begin{lemma}\label{th measurability}
Let a process $Z$ be adapted to a non increasing filtration  $\mathcal{G}=(\mathcal{G}_t,t\geq0)$. Let $\Gamma$ be an optional time with respect to $\mathcal{G}$, i.e. for all $t\geq0$, the event $\{\Gamma> t\}\in\mathcal{G}_t$. The following holds.
\begin{enumerate}
\item The set $\mathcal{G}_{\Gamma}:=\{A\in\mathcal{G}_0\:: A\cap\{\Gamma> t\}\in\mathcal{G}_t\}$ is a sigma-algebra.
\item 
If all the paths of $Z$ are either continuous on the right  or on the left, then the r.v. $Z_{\Gamma}\1_{\Gamma <\infty}$ is $\mathcal{G}_{\Gamma}$ measurable.
\item Suppose that $\mathcal{G}$ is continuous on the right; that is, for all $t$, 
$\mathcal{G}_t=\sigma\left(\underset{s>t}{\cup}\mathcal{G}_s\right)$. If $Z$ is a backward martingale with respect to $\mathcal{G}$, then for all $t$, the right hand and left hand limits $Z_{t^+}$, $Z_{t^-}$ exist a.s. 
Moreover, the process $t\mapsto Z_{(\Gamma\vee t)^+}-\Delta_{\Gamma}\1_{\Gamma> t}$ is a backward martingale with respect to the completed filtration
$\overline{\mathcal{G}}$, with $\Delta_{\Gamma}=Z_{\Gamma^+}-Z_{\Gamma^-}$.
\end{enumerate}
\end{lemma}

In the latter, the completed filtration is defined by 
$\overline{\mathcal{G}}_t=\sigma(\mathcal{N}\cup\mathcal{G}_t)$ and $\mathcal{N}=\{A\in\mathcal{G}_0\:: P(A)=0\}$.
 The following is also a useful tool.

\begin{lemma}\label{th stopping time}
\begin{enumerate}
\item If a process $Z$ is such that $Z_{s+t}=\phi_{s,t}(Z_{s})$ for all $t,s\geq0$, then $\tau(Z_0)=:\Gamma$ is an optional time with respect to the natural non increasing filtration $\mathcal{F}^{Z}$ of $Z$. 
Moreover, $\mathcal{F}^{Z}_0=\mathcal{F}^{Z}_{\Gamma}$.

\item Suppose that $\{\Gamma\leq t\}\in\mathcal{F}^Z\cap\mathcal{F}^{Z'}$ for some $t\geq0$. If $Z'_t\1_{\Gamma\leq t}=Z_t\1_{\Gamma\leq t}$, then 
$\E[F|Z'_t]\1_{\Gamma\leq t}=\E[F|Z_t]\1_{\Gamma\leq t}$ for all integrable r.v. $F$.
\end{enumerate}
\end{lemma}
From proposition \ref{th turbulent flows}, this lemma is satisfied by $Z^1,Z^2,Z^3,Z^4$ and $X$, taking the role of $Z,Z'$.\\

\textbf{Proof.} We begin with the first assertion. $u^-(\cdot,t),u^+(\cdot,t)$ are Borel functions and it is well known that if $u$ is discontinuous in $(Z_t,t)$, it is also discontinuous in $(Z_{t+s},t+s)$. Then, 
$$\{\Gamma\leq t\}=\{u^-(Z_t,t)\neq u^+(Z_t,t)\}\cup[\{u^-(Z_t,t)= u^+(Z_t,t)\}\cap\{\Gamma= t\}]\;.$$
Since
\begin{align*}
\{u^-(Z_t,t)= u^+(Z_t,t)\}&\cap\{\Gamma= t\}=\{u^-(Z_t,t)= u^+(Z_t,t)\}\cap\\
&\left[\underset{n\geq1}{\cap}\{u^-(Z_{t+1/n},t+1/n)\neq u^+(Z_{t+1/n},t+1/n)\}\right],
\end{align*}
the proof of the first assertion is done.\\
Remark that $Z_{t+1/n}=\phi_{t,1/n}(Z_t)$. So $\{\Gamma\leq t\}=Z_t^{-1}(A_t)$, with
\begin{align*}
A_t=\{u^-(\cdot,t)\neq u^+(\cdot,t)\}\cup&\bigg(\{u^-(\cdot,t)= u^+(\cdot,t)\}\cap\\
&\left[\underset{n\geq1}{\cap}\{u^-(\phi_{t,1/n},t+1/n)\neq u^+(\phi_{t,1/n},t+1/n)\}\right]\bigg)
\end{align*}

Now we show that $\mathcal{F}^{Z}_0=\mathcal{F}^{Z}_{\Gamma}$. First remark that if 
$\{b\}\neq\phi_{0,t}^{-1}(\phi_{0,t}(b))$, then $\tau(b)\leq t$. Thus for all Borel subset $B$ and $t\geq0$, we have
$B\cap\{\tau>t\}=\phi_{0,t}^{-1}(\phi_{0,t}(B))\cap\{\tau>t\}$ and
\begin{align*}
Z_0^{-1}(B)\cap\{\tau(Z_0)>t\}&=Z_t^{-1}(\phi_{0,t}(B))\cap\{\tau(Z_0)>t\}\;,\\
Z_0^{-1}(B)\cap\{\Gamma>t\}&=Z_t^{-1}(\phi_{0,t}(B))\cap\{\Gamma>t\}\in\mathcal{F}^{Z}_t\;.
\end{align*}
This means that $Z_0^{-1}(B)\in\mathcal{F}^{Z}_{\Gamma}$.\\

For the second assertion, since $Z_t\1_{\gamma\leq t}=Z'_t\1_{\Gamma\leq t}$, it is easy to see that $\E[F|Z'_t]\1_{\Gamma\leq t}$ is 
$\sigma(Z'_t)\cap\sigma(Z_t)$ measurable; for all bounded Borel function $h$, 
\begin{align*}
\E\big(h(Z_t)\E[F|Z'_t]\1_{\Gamma\leq t}\big)=\E\big(h(Z'_t)\E[F|Z'_t]\1_{\Gamma\leq t}\big)
=\E\big(h(Z'_t)F\1_{\Gamma\leq t}\big)\\
=\E\big(h(Z_t)F\1_{\Gamma\leq t}\big)=\E\big(h(Z_t)\E[F|Z_t]\1_{\Gamma\leq t}\big).
\end{align*}
Hence, $\E[F|Z'_t]\1_{\Gamma\leq t}=\E[F|Z_t]\1_{\Gamma\leq t}$ a.s.

\begin{corollary}[Turbulence semi-martingales]\label{th turbulence semimartingales}
Let $Z$ stand independently for $Z^1,Z^2,Z^3$ or $Z^4$. The process $Z$ is Markovian. Suppose that the support of $\mu_0$ is $\mathcal{S}=\R$.
\begin{enumerate}
\item $t\mapsto u(Z_{t},t)\1_{\gamma>t}$ is a bounded variational process adapted to the natural non increasing filtration $\mathcal{F}^Z$ of $Z$.
\item $t\mapsto u(Z_{t},t)\1_{\gamma\leq t}$ is a backward \textit{c\`adl\`ag} semi-martingale with respect to $\mathcal{F}^Z$. Moreover,
$t\mapsto u(Z_{t},t)\1_{\gamma\leq t}+\E[u_0(X_0)|Z_0]\1_{t<\gamma}$ is a backward martingale.
\item $t\mapsto u(Z_{t},t)$ is a backward \textit{c\`adl\`ag} semi-martingale with respect to $\mathcal{F}^Z$. Moreover, 
$t\mapsto u(Z_{t},t)-\left(u_0(Z_0)-\E[u_0(X_0)|Z_0]\right)\1_{t<\gamma}$ is a backward martingale.
\end{enumerate}
\end{corollary}

\noindent\textbf{Proof.}  \textbf{1)} Obviously, $u(Z_{t},t)\1_{\gamma> t}=u_0(Z_0)\1_{\gamma> t}$. \\ 
\noindent\textbf{2)} For all $t$, 
\begin{align*}
u(Z_{t},t)\1_{\gamma\leq t}&=u(X_{t},t)\1_{\gamma\leq t}=\E[u_0(X_0)|X_t]\1_{\gamma\leq t}
=\underbrace{\E[u_0(X_0)|Z_t]}_{M_t}\1_{\gamma\leq t}=M_{t^+}\1_{\gamma\leq t}\\
&=M_{(\gamma\vee t)^+}-M_{\gamma^+}\1_{t<\gamma}\;.
\end{align*}
But the process $t\mapsto M_{\gamma^+}\1_{t<\gamma}$ is adapted to $\mathcal{F}^Z$. 
Then, according to lemma \ref{th measurability}, the process $t\mapsto u(Z_{t},t)\1_{\gamma\leq t}+ M_{\gamma^-}\1_{t<\gamma}$ is a backward martingale with respect to $\mathcal{F}^Z$.
Moreover 
\begin{align*}
M_{\gamma^-}&=\underset{\underset{s<\gamma}{s\to\gamma}}{\lim}M_s\1_{s<\gamma}
=\underset{\underset{s<\gamma}{s\to\gamma}}{\lim}\E\left[M_0\1_{s<\gamma}|Z_s\right]\\
&=\underset{\underset{s<\gamma}{s\to\gamma}}{\lim}M_0\1_{s<\gamma}=M_0=\E[u_0(X_0)|Z_0]\;.
\end{align*}
\noindent\textbf{3)} $u(Z_{t},t)-(u_0(Z_0)-\E[u_0(X_0)|Z_0])\1_{t<\gamma}=
M_{(\gamma\vee t)^+}-[M_{\gamma^+}-M_{\gamma^-}]\1_{t<\gamma}$. \\

Now, in order to analyse the process of \cite{spontaneous-Drivas-Eyink}, let us define, 
for all $(a,t,s,r)\in\mathbb{R}\times \mathbb{R}^{3}_{+}$
\begin{align*}
f(a,t,r) & =  \left\{  \begin{array}{ll}
\phi_{0,t} (a) &  \textrm{if  $t\geq r$}\\
\phi_{0,r} (a)-(r-t) u^{-}\big(\phi_{0,r} (a),\,r\big)  &  \textrm{if\,$t<r$}\\
\end{array}  \right.\\
g(a,t,r) & =  \left\{  \begin{array}{ll}
\phi_{0,t} (a) &  \textrm{if  $t\geq r$}\\
\phi_{0,r} (a)-(r-t) u^{+}\big(\phi_{0,r} (a),\,r\big)  &  \textrm{if\,$t<r$}\\
\end{array}  \right.
\end{align*} 
The definition of \cite{spontaneous-Drivas-Eyink} is the following : 
\begin{eqnarray*}\label{process def}
Y_{t}(a)  =  \left\{  
\begin{array}{ll}
 f(a,t,\tau(a)) &  \textrm{if $a$ enters in the shock from the left} \\
 g(a,t,\tau(a)) &  \textrm{if $a$ enters in the shock from the right}. 
\end{array}  \right.
\end{eqnarray*}
In this definition, if a turbulent interval $[\alpha,\beta]$ is also a cluster at time $\tau(\alpha)$, then the whole interval $[\alpha,\beta[$ is considered as entering in the shock from the left, and $\beta$ is considered as entering in the shock from the right.

This definition is however ambiguous since it supposes that a particle can hurt only one discontinuity line at a time. In the sequel, we call this \textit{the assumption of simple shocks}. 
In order to remove the ambiguity, we suggest to slightly modify the definition :
\begin{eqnarray*}
Y_{t}(a)  =  \left\{  
\begin{array}{ll}
 f(a,t,\tau(a)) &  \textrm{if $a$ enters in the shock from the left} \\
								&   \textrm{of all concerned discontinuity lines}\\
 g(a,t,\tau(a)) &  \textrm{otherwise.} 
\end{array}  \right.
\end{eqnarray*}
Of course, this definition matches the underlying \textit{assumption of simple shocks} of \cite{spontaneous-Drivas-Eyink}. In that paper, the authors considered Lebesgue measure as initial distribution $\mu_0$ of the matter. 

\begin{figure}
\centering
\includegraphics[width=1.25\linewidth]{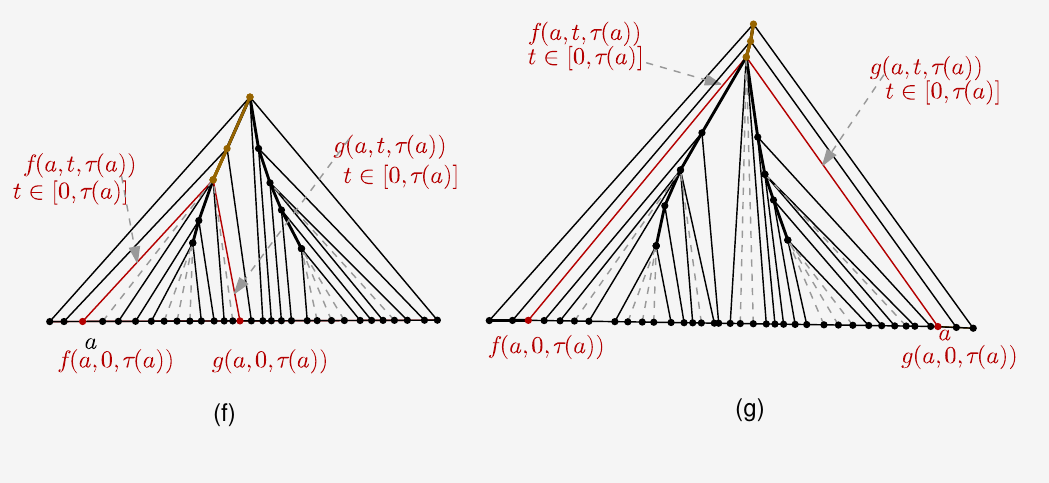}
\caption[trajec delta shock]{\scriptsize   {Delta-shock.}}
\label{Delta-shock}
\end{figure}

\begin{proposition}[Link with delta-shocks]\label{th delta-prime-flow}
Suppose that the support of $\mu_0$ is $\mathcal{S}=\R$.
For all $0\leq t\leq r$ and $a\in\R$, let $x=\phi_{0,r}(a)$. We have
\begin{align}
\label{shock entering left}
f(a,t,r)  &=\phi_{0,t}(\alpha(x,0,r))\\
&=\alpha(x,0,r)+tu_0(\alpha(x,0,r))\;.\nonumber\\
\label{shock entering right}
g(a,t,r) & =\phi_{0,t}(\beta(x,0,r))\\
&=\beta(x,0,r)+tu_0(\beta(x,0,r))\;.\nonumber
\end{align} 
If moreover $\alpha(x,0,r)<\beta(x,0,r)$, then $$r=\tau(\alpha(x,0,r))=\tau(\beta(x,0,r))=\tau(a)\;.$$
\end{proposition}
Thus, $[0,\tau(a)]\ni t\mapsto f(a,t,\tau(a)),\,g(a,t,\tau(a))$ draw a delta-shock (figure \ref{Delta-shock}). 

We can (more generally) consider $Y$ on abstract set $\Omega$ (instead of $\R$) : 
\begin{eqnarray}
Y_{t}  =  \left\{  
\begin{array}{ll}
 f(X_0,t,\gamma) &  \textrm{if $X_0$ enters in the shock from the left} \\
								&   \textrm{of all concerned discontinuity lines}\\
 g(X_0,t,\gamma) &  \textrm{otherwise.} 
\end{array}  \right.
\end{eqnarray}
The process $Y$ then chooses one segment of delta-shock, at random. Note that the event
"$X_0$ enters the shock from the left" coincides with\\ 
$\{Z^1_0=Z^3_0\}\cap[\{Z^2_0=Z^4_0,\,X_0\neq Z^2_0\}\cup\{Z^2_0\neq Z^4_0\}]$.
\begin{proposition}[Motion on delta-shock]\label{th Motion on one wing}
Suppose that the support of $\mu_0$ is $\mathcal{S}=\R$. 
\begin{enumerate}
\item $\tau(Y_0)=\tau(X_0)=\gamma$ and 
$Y_t=\left\{
\begin{array}{ll}
 Y_0+tu_0(Y_0)&  \textrm{if } t\leq\gamma\\
Z_t=X_t&  \textrm{if } t\geq\gamma\;. 
\end{array}\right.
$
\begin{align}\label{proc papillon}
\forall\,s\geq0,\,\forall\,t\geq0\;,\quad  Y_{s+t}=\phi_{s,t}(Y_s)=\left\{  
\begin{array}{ll}
 Z^3_{s+t}&  \textrm{if } Z^1_0=Z^3_0\\
Z^4_{s+t}&  \textrm{if } Z^1_0\neq Z^3_0\;.
\end{array}\right.
\end{align} 
\item If the assumption of simple shocks holds, then two events occur :
\begin{enumerate}
\item $Y_0=Z^3_0=Z^1_0$ and in this case, for all $t\geq0$, $Y_t=Z^3_t=Z^1_t$.
\item $Y_0=Z^4_0=Z^2_0$ and in this case, for all $t\geq0$, $Y_t=Z^4_t=Z^2_t$.
\end{enumerate}
\item If $\gamma$ has no atom, i.e. $P(\gamma=t)=0$ for all $t$, then $X\equiv Z^1\equiv Z^2$. And in this case, all the turbulent intervals are reduced to single points. If moreover the assumption of simple shocks holds, then $Y\equiv X\equiv Z^1\equiv Z^2$.
\end{enumerate}
\end{proposition}

\begin{remark}[Delta-shock velocity as semi-martingale]\label{rem random wing semimartingale}
\begin{enumerate}
\item From the first assertion of proposition \ref{th Motion on one wing}, the process $t\mapsto u(Y_{t},t)$ satisfies corollary 
\ref{th turbulence semimartingales} (with $Z=Y$) and is a backward semi-martingale of $\mathcal{F}^Y$.
\item The last assertion of proposition \ref{th Motion on one wing} shows that the martingale $t\mapsto u(Y_{t},t)$ of \cite{spontaneous-Drivas-Eyink} (which required the assumption of uniformity of the law of $\gamma$) was in fact already obtained in \cite{MoutsingaPhD}. The case of \cite{spontaneous-Drivas-Eyink} is a particularity where all the turbulent intervals are reduced to single points. 
\end{enumerate}
\end{remark}

Now we precise, under more general assumptions, when the velocity of turbulence is a martingale.

\subsection{Martingales and undetectability of turbulence}
In this part, we show that the martingality of the velocity of turbulence implies that any turbulent interval is a single point (\textit{single turbulent point}). 

\begin{corollary}[Turbulence martingales and prime-delta-shocks]\label{Th martingales}
Let $Z$ stand independently for $Z^1,Z^2,Z^3,Z^4$, or $Y$. Suppose that the support of $\mu_0$ is $\mathcal{S}=\R$. 

The process $t\mapsto u(Z_{t},t)$  is a martingale if and only if a.e. 
$$X\equiv Z\equiv Z^1\equiv Z^2.$$
\end{corollary}

\textbf{Proof.}  The semi-martingale is a martingale iff its bounded variational part
$t\mapsto \left(u_0(Z_0)-\E[u_0(X_0)|Z_0]\right)\1_{t<\gamma}$  is constant. Letting $t$ tend to 0 and $+\infty$ respectively gives
\begin{align*}
\left(u_0(Z_0)-\E[u_0(X_0)|Z_0]\right)\1_{\gamma>0}
&=\left(u_0(Z_0)-\E[u_0(X_0)|Z_0]\right)\1_{\gamma=\infty}=0
\end{align*}
since $\gamma$ is $\sigma(Z_0)$-measurable and $X_0\1_{\gamma=\infty}=Z_0\1_{\gamma=\infty}$. 
(In the same way, $X_0\1_{\gamma=0}=Z_0\1_{\gamma=0}$.) %; so we get $u_0(Z_0)-\E[u_0(X_0)|Z_0]$
Using the fact that $X_0+\gamma u_0(X_0)=Z_0+\gamma u_0(Z_0)$, the NSC to have a martingale becomes
$$\E\left[\gamma^{-1}(X_0-Z_0)\1_{\gamma>0}|Z_0\right]=\E\left[(u_0(Z_0)-u_0(X_0))\1_{\gamma>0}|Z_0\right]=0.$$ 
Let us now study each case of $Z$. 
\begin{enumerate}
\item For $Z=Z^3$ : we have 
$\E[\gamma^{-1}(X_0-Z^3_0)\1_{\gamma>0}].$ 
But $X_0\geq Z^3_0$. So $X_0=Z^3_0=Z^1_0$ a.e.
Thus $X\equiv Z^3\equiv Z^1$ a.e.\\
Now, we show that $Z^1_0=Z^2_0$. 
If $[Z^1_0,Z^2_0]=[\alpha,\beta]$, then $\mu_0(]\alpha,\beta])=P(\alpha <X_0 \leq\beta)\leq P(X_0\neq Z^1_0)=0$. This implies that $\alpha=\beta$ since there is no vacuum in the support.  Hence $Z^1_0=Z^2_0$. 

\item For $Z=Z^4$, resp. $Z=Z^1$, resp. $Z=Z^2$  : Analogous to previous case.

\item For $Z=Y$ : let $L$ be the set of particles which enter the shock on the left. We have 
$\E[\gamma^{-1}(X_0-Y_0)\1_{\gamma>0}\1_L(Y_0)]=0$, with  $(X_0-Y_0)\1_L(Y_0)\geq0$, so $(X_0-Y_0)\1_L(Y_0)=0$. In the same way, 
$(X_0-Y_0)\1_{L^c}(Y_0)=0$ and we get $X_0=Y_0$. \\
Now, we show that $Z^1_0=Z^2_0$. If $[Z^1_0,Z^2_0]=[\alpha,\beta]$, then $\mu_0(]\alpha,\beta[)=P(\alpha <X_0 <\beta)\leq P(X_0\neq Y_0)=0$. This implies that $\alpha=\beta$.  Hence $Z^1_0=Z^2_0$. 
\end{enumerate}

% % % % % % % % % % % % % % % % % % % % % % % % % % % % % % %
 {\small
  \bibliographystyle{abbrv} 	
  \bibliography{biblio_Burgers} 	 
   	  	  	   	  	  	 }
\end{document}